\documentclass[11pt]{amsart}
\usepackage{amssymb,amsmath,amsfonts,amscd,euscript}

\newcommand{\nc}{\newcommand}

\numberwithin{equation}{section}
\newtheorem{thm}{Theorem}[section]
\newtheorem{prop}[thm]{Proposition}
\newtheorem{lem}[thm]{Lemma}
\newtheorem{cor}[thm]{Corollary}
\theoremstyle{remark}
\newtheorem{rem}[thm]{Remark}

\newtheorem{example}[thm]{Example}
\newtheorem{dfn}[thm]{Definition}

\nc{\gl}{\mathfrak{gl}}
\nc{\GL}{\mathfrak{GL}}
\nc{\g}{\mathfrak{g}}
\nc{\gh}{\widehat\g}
\nc{\h}{\mathfrak{h}}
\nc{\la}{\lambda}
\nc{\al}{\alpha }
\nc{\be}{\beta }
\nc{\ve}{\varepsilon }
\nc{\om}{\omega }

\nc{\ta}{\theta}
\nc{\veps}{\varepsilon}
\nc{\ch}{{\mathop {\rm ch}}}
\nc{\Tr}{{\mathop {\rm Tr}\,}}
\nc{\Id}{{\mathop {\rm Id}}}
\nc{\ad}{{\mathop {\rm ad}}}
\nc{\bra}{\langle}
\nc{\ket}{\rangle}
\nc{\x}{{\bf x}}
\nc{\bs}{{\bf s}}
\nc{\bp}{{\bf p}}
\nc{\bc}{{\bf c}}
\nc{\pa}{\partial}
\nc{\ld}{\ldots}
\nc{\cd}{\cdots}
\nc{\hk}{\hookrightarrow}
\nc{\T}{\otimes}
\newcommand{\bea}{\begin{equation}}
\newcommand{\ena}{\end{equation}}
\nc{\gr}{\mathrm{gr}}
\nc{\ov}{\overline}

\nc{\cO}{\mathcal O}
\nc{\msl}{\mathfrak{sl}}
\nc{\mgl}{\mathfrak{gl}}
\nc{\U}{\mathrm U}
\nc{\V}{\EuScript V}
\nc{\bH}{\EuScript H}
\nc{\Res}{\mathrm{Res\ }}

\newcommand{\bC}{{\mathbb C}}
\newcommand{\bZ}{{\mathbb Z}}

\newcommand{\bP}{{\mathbb P}}
\newcommand{\bG}{{\mathbb G}}

\newcommand{\fb}{{\mathfrak b}}

\newcommand{\fn}{{\mathfrak n}}

\newcommand{\fu}{{\mathfrak u}}

\newcommand{\Fl}{\EuScript{F}}

\newcommand{\spa}{\mathrm{span}}

\begin{document}

\title[Degenerate  flag varieties]
{Degenerate flag varieties and the median Genocchi numbers}

\author{Evgeny Feigin}
\address{Evgeny Feigin:\newline
Mathematical Department, University Hisger School of Economics,\newline
20 Myasnitskaya st, 101000, Moscow, Russia\newline
{\it and }\newline
Tamm Theory Division,
Lebedev Physics Institute,\newline
Leninisky prospect, 53,
119991, Moscow, Russia
}
\email{evgfeig@gmail.com}

\begin{abstract}
We study the $\bG_a^M$ degenerations $\Fl^a_\la$ of the type $A$ flag varieties $\Fl_\la$.
We describe these degenerations explicitly as subvarieties in the products of Grassmanians.
We construct cell decompositions of $\Fl^a_\la$ and show that for complete flags the number of cells
is equal to the normalized median Genocchi numbers $h_n$.
This leads to a new combinatorial definition of the  numbers $h_n$.
We also compute the Poincar\' e polynomials of the complete degenerate flag varieties
via a natural statistics
on the set of Dellac's configurations, similar to the length statistics on the
set of permutations. We thus obtain a natural $q$-version of the normalized median Genocchi numbers.
\end{abstract}

\maketitle

\section*{Introduction}
Let $\g=\msl_n$, $G=SL_n$. Fix the Cartan decomposition $\g=\fb\oplus \fn^-$, where $\fb$
is a Borel subalgebra, $\fb=\fn\oplus\h$.
In \cite{Fe3} we considered the degenerate algebra $\g^a=\fb\oplus (\fn^-)^a$, where $(\fn^-)^a$
is an abelian Lie algebra isomorphic to $\fn^-$ as a vector space. The corresponding Lie group
is a semi-direct product $G^a=B\ltimes \bG_a^M$, where $\bG_a$ is the additive group of the field and
$M=\dim \fn$. For a dominant integral weight $\la$ let $V_\la$ be the highest weight $\la$
irreducible $\g$-module with a highest weight vector $v_\la$.
The increasing PBW filtration $F_\bullet$ on $V_\la$ is defined as follows:
\[
F_0=\bC v_\la,\ F_{s+1}=\mathrm{span}\{xv:\ x\in\g, v\in F_s \}, s\ge 0
\]
(see \cite{Fe1}, \cite{Fe2}, \cite{FFoL1}, \cite{FFoL2}, \cite{Kum2}).
The associated graded space $V_\la^a=F_0\oplus F_1/F_0\oplus F_2/F_1\oplus\dots$ can be naturally
endowed with the structure of a $\g^a$- and $G^a$-module.
A degenerate flag variety
$\Fl^a_\la$ is a subvariety in $\bP(V_\la^a)$ defined by $\Fl^a_\la=\overline{\bG_a^M\cdot \bC v_\la}$.
These are the $\bG_a^M$-degenerations of the classical (generalized) flag varieties $\Fl_\la$
(see \cite{A}, \cite{AS}, \cite{Fe3}, \cite{HT}).
For example, $\Fl^a_{\om_d}\simeq Gr(d,n)$
for all fundamental weights. Recall also that in the classical case (for $\g=\msl_n$) the varieties
$\Fl_\la=G\cdot \bC v_\la\hk \bP(V_\la)$ are the usual flag varieties (maybe partial).
In particular, if $\la$ is regular, i.e. $(\la,\om_d)>0$ for all $d$, then $\Fl_\la$ is isomorphic
to the variety $\Fl_n$ of complete flags in $n$-dimensional space $V$. Fix a basis  $v_1,\dots,v_n$
of $V$.

For all weights $\la$, $\mu$ there exists an embedding of $G^a$-modules
$V_{\la+\mu}^a\hk V_{\la}^a\T V_{\mu}^a$ sending $v_{\la+\mu}$ to $v_\la\T v_\mu$ (see \cite{FFoL1}, \cite{FFoL2}).  \
This induces the embedding
of varieties $\Fl_{\la+\mu}^a\hk \Fl_{\la}^a\times \Fl_{\mu}^a$. Thus for any $\la$ we obtain an embedding
of $\Fl^a_\la$ into the product of Grassmanians. Our first result is an explicit description of this embedding.
We state the theorem here for complete flag varieties $\Fl^a_n$.
For this we need one more piece of notations. Let $pr_d:V\to V$ be the projection along
the space $\bC v_d$ to the linear span of the vectors $v_i$, $i\ne d$.

\begin{thm}\label{intro1}
The image of the embedding of  the variety $\Fl^a_n$ in the product $\prod_{d=1}^{n-1} Gr(d,n)$ is equal to the set
of chains of subspaces $(V_1,\dots,V_{n-1})$,  $V_d\in Gr(d,n)$ such that
\[
pr_{d+1}(V_d)\hk V_{d+1}, \quad 1\le d\le n-2.
\]
\end{thm}

Our next goal is to compute the Poincar\' e polynomial of $\Fl_n^a$. Recall that in the classical case
the flag variety $\Fl_n$ can be written as a disjoint union of $n!$ cells, each cell being
associated with a torus fixed point. The fixed points are labeled by permutations from $S_n$.
The length statistics $\sigma\to l(\sigma)$ gives the complex dimension of the cells.
Therefore, the  Poincar\' e polynomial $P_{\Fl_n}(t)$ of $\Fl_n$ is equal to
$P_{\Fl_n}(t)=\sum_{\sigma\in S_n} t^{2l(\sigma)}$.

As an immediate corollary of Theorem \ref{intro1} we obtain that the fixed points of the torus
$T\subset G^a$ action on $\Fl^a_n$
are labeled by the sequences $I^1,\dots, I^{n-1}$, $I^d\subset \{1,\dots,n\}$, $\# I^d=d$,
satisfying
\begin{equation}\label{introgen}
I^d\setminus\{d+1\}\hk I^{d+1},\qquad  d=1,\dots,n-2.
\end{equation}
(Note that this set of sequences
has a subset with $I^d\hk I^{d+1}$, which can be naturally identified with the permutations $S_n$).
Our first task is to compute the number of such fixed points. To this end, recall the normalized median
Genocchi numbers $h_n$, $n=1,2,\dots$ (also referred to as the normalized Genocchi numbers of second kind).
These numbers have several definitions \cite{De}, \cite{Du}, \cite{DR}, \cite{DZ}, \cite{G}, \cite{Kr}, \cite{Vien}
(see section \ref{sec3} for a review). Here we give
the Dellac definition, which is the earliest one and
which fits our construction in the best way.

Consider a rectangle with $n$ columns and
$2n$ rows. It contains $n\times 2n$ boxes labeled by pairs $(l,j)$,
with $l=1,\dots,n$ and $j=1,\dots,2n$. A Dellac  configuration $D$ is a subset of boxes,
subject to the following conditions: first, each column contains exactly two boxes from $D$ and
each row contains exactly one box from $D$, and, second, if the $(l,j)$-th box is in $D$,
then $l\le j\le n+l$.
Let $DC_n$ be the set of such configurations. Then $h_n$ is the number of elements in $DC_n$.
The first several median Genocchi numbers (starting from $h_1$) are as follows: $1,2,7,38,295$.
For instance, the two Dellac configurations for $n=2$ are as follows:
(we specify boxes in a configuration by putting  fat dots inside)
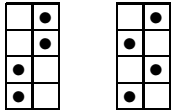
\[
\begin{picture}(20,40)
\put(0,0){\line(1,0){20}}
\put(0,10){\line(1,0){20}}
\put(0,20){\line(1,0){20}}
\put(0,30){\line(1,0){20}}
\put(0,40){\line(1,0){20}}

\put(0,0){\line(0,1){40}}
\put(10,0){\line(0,1){40}}
\put(20,0){\line(0,1){40}}

\put(2,2){$\bullet$}
\put(2,12){$\bullet$}
\put(12,22){$\bullet$}
\put(12,32){$\bullet$}
\end{picture}\quad\quad
\begin{picture}(20,40)
\put(0,0){\line(1,0){20}}
\put(0,10){\line(1,0){20}}
\put(0,20){\line(1,0){20}}
\put(0,30){\line(1,0){20}}
\put(0,40){\line(1,0){20}}

\put(0,0){\line(0,1){40}}
\put(10,0){\line(0,1){40}}
\put(20,0){\line(0,1){40}}

\put(2,2){$\bullet$}
\put(2,22){$\bullet$}
\put(12,12){$\bullet$}
\put(12,32){$\bullet$}
\end{picture}
\]

We prove the following theorem:
\begin{thm}
The number of sequences $I^1,\dots,I^{n-1}$ as above, satisfying \eqref{introgen} is equal to $h_n$.
\end{thm}
We also prove that the Dellac definition \cite{De} is equivalent to the Dumont-Kreweras definition
\cite{Du}, \cite{Kr}
(this fact is known to experts \cite{G},\cite{S} but we were unable to find the proof in the literature).

Recall that the length of a permutation $\sigma\in S_n$ can be defined as the number of pairs
$1\le l_1<l_2\le n$ satisfying $\sigma(l_1)>\sigma(l_2)$. We define a length $l(D)$ of a Dellac configuration
$D$ as the number of squares $(l_1,j_1), (l_2,j_2)\in D$ such that $l_1<l_2$ and $j_1>j_2$.
We prove the following theorem:
\begin{thm}
The Poincar\'e polynomial $P_{\Fl_n^a}(t)$ is given by $\sum_{D\in DC_n} t^{2l(D)}$.
\end{thm}

Our paper is organized in the following way:\\
In Section \ref{sec1} we recall main definitions and theorems from \cite{Fe3},\\
In Section \ref{sec2} we describe explicitly the image of the embedding of the varieties $\Fl^a_\la$
into the product of Grassmanians and construct the cell decomposition of $\Fl^a_\la$.\\
In Section \ref{sec3} we study the combinatorics of the median Genocchi numbers and compute
the Poincar\'e polynomials of the complete degenerate flag varieties.

\section{PBW deformation}\label{sec1}
\subsection{Definitions}
We first recall basic definitions and constructions from \cite{FFoL1} and \cite{Fe3}.
Let $\g$ be a simple Lie algebra with the Cartan decomposition
$\g=\fn\oplus\h\oplus \fn^-$. We denote by $M$ the number of positive roots of $\g$, i.e.
$M=\dim\fn$. Let $\fb=\fn\oplus\h$ be a Borel subalgebra.
Then the
deformed algebra $\g^a$ is defined as a sum of two subalgebras $\g^a=\fb\oplus (\fn^-)^a$,
where $(\fn^-)^a$ is an abelian Lie algebra isomorphic to $\fn^-$ as a vector space.
The subalgebra $(\fn^-)^a\hk \g^a$ is an abelian ideal and the action of
$\fb$ on $(\fn^-)^a$ is induced from the identification $(\fn^-)^a\simeq \g/\fb$.

Let $G$ be the Lie group of the Lie algebra $\g$. Let $N, T, N^-, B$ be the Lie groups of
the Lie algebras $\fn$, $\h$, $\fn^-$, $\fb$.
The deformed Lie group $G^a$ is defined as a semi-direct product of $B$ and
the normal subgroup $\bG_a^M$,
where $\bG_a$ is the additive group of the field (thus $\bG_a^M$ is the Lie group
of the Lie algebra $(\fn^-)^a$). The Borel group $B$ acts on the vector space
$(\fn^-)^a\simeq \g/\fb$ via the restriction of the adjoint action and therefore there
exists a natural homomorphism from $B$ to $Aut(\bG_a^M)$, defining the semi-direct product
$G^a=B\ltimes \bG_a^M$.

For a dominant integral weight $\la$ we denote by $V_\la$ the corresponding irreducible highest weight
$\g$-module with a highest weight vector $v_\la$.
The Lie algebra $\g^a$ and the Lie group $G^a$ act on the deformed representations $V_\la^a$,
where $\la$ are dominant integral weights of $\g$. The representations $V_\la^a$
are defined as associated graded $gr_\bullet V_\la$ of the representation $V_\la$ with
respect to the PBW filtration $F_s$:
\[
F_s=\mathrm{span}\{x_1\dots x_lv_\la:\ x_i\in\g, l\le s\}.
\]
So $V_\la^a=\bigoplus_{s\ge 0} V_\la^a(s)$, where $V_\la^a(0)=\bC v_\la$ and
$V_\la^a(s)=F_s/F_{s-1}$ for $s>0$. It is easy to see that the action of $\fn^-$
on $V_\la$ becomes abelian on $V_\la^a$ (i.e. it induces the action of $(\fn^-)^a$)
and the action of the Borel subalgebra induces the action of (the same algebra) $\fb$.
The actions of $(\fn^-)^a$ and $\fb$ glue together to the action of $\g^a$.\
\begin{rem}
Let $\tilde\g^a=\g^a\oplus \bC p$ be the central, extension of $\g^a$ with a single element
$p$ subject to the relations $[p,\fb]=0$, $[p,f_\al]=f_\al$ for any positive root $\al$ and the corresponding
weight element $f_\al\in(\fn^-)^a$. Thus the Cartan subalgebra of $\tilde\g^a$ has one extra
dimension. We note that the $\g^a$-module structure of $V_\la^a$ naturally lifts to the structure
of representation of $\tilde\g^a$ by setting $pv_\la=0$ (in general, $p|_{V^a_\la(s)}=s$).
An eigenvalue of the operator $p$ is sometimes referred to as a PBW degree.
The character of $V_\la^a$ with respect to $\h\oplus \bC p$ was computed in \cite{FFoL1} for $\msl_n$
and in \cite{FFoL2} for symplectic Lie algebras.
We denote the
Lie group of $\tilde\g^a$ by $\tilde G^a$, which differs from $G^a$ by an additional $\bC^*$.
\end{rem}

Consider the action of $G^a$ on the projective
space $\bP(V_\la^a)$.
Recall that in the classical situation the (generalized) flag varieties are defined as
$\Fl_\la=G\cdot \bC v_\la\hk \bP(V_\la)$ (see \cite{Kum1}).
The degenerate flag varieties $\Fl^a_\la\hk \bP(V_\la^a)$ are defined as the closures
of the $G^a$ orbit (or, equivalently, of the $\bG_a^M$ orbit) of the line $\bC v_\la$.
We note that in the classical case the orbit $G\cdot \bC v_\la$ already covers the whole flag variety.
This is not true in the degenerate case: the orbit $G^a  \cdot \bC v_\la$ is an affine cell, whose
closure gives a projective singular variety $\Fl^a_\la$.

\subsection{The type $A$ case}
From now on we assume that $\g=\msl_n$ and $G=SL_n$.
Then all positive roots are of the form
\[
\al_{i,j}=\al_i+\dots +\al_j,\ 1\le i\le j\le n-1
\]
(for instance, $\al_{i,i}=\al_i$ are the simple roots).
We denote by $f_{i,j}=f_{\al_{i,j}}\in \fn^-$ and $e_{i,j}=e_{\al_{i,j}}\in \fn$
the corresponding root elements.
We have
$\Fl^a_{\om_d}\simeq \Fl_{\om_d}\simeq Gr(d,n)$. The reason why the degenerate flag
varieties are isomorphic to the non-degenerate ones for fundamental weights is that
the radicals in $\msl_n$, corresponding to $\om_d$, are abelian. In other words, define the set
of positive roots
\[
R_d=\{\al_{i,j}:\ 1\le i\le d\le j\le n-1\}.
\]
Define the subalgebra $\fu^-_d=\mathrm{span}\{f_\al:\ \al\in R_d\}$.
Then $\fu^-_d$ is abelian and $V_{\om_d}=U(\fu^-_d)\cdot v_\la$.

\begin{rem}
Let us explain the difference between the structure of $\g$-module on $V_{\om_d}$ and
the structure of $\g^a$-module on $V_{\om_d}^a$.  The operators $f_\al$ act trivially on $V_{\om_d}^a$
unless $\al\in R_d$. Also, $e_\al$ act trivially on $V_{\om_d}^a$ if $\al\in R_d$.
Therefore, $\g^a$ acts on $V_{\om_d}^a$ via the projection to the subalgebra
\begin{equation}\label{gd}
\g^a_d=\fu^-_d\oplus \h\oplus \mathrm{span}\{e_\al:\ \al\notin R_d\}.
\end{equation}
Similarly, the group $G^a$ acts on $Gr(d,n)$ via the surjection to the Lie group
of $\g^a_d$. In particular, the group $G^a$ does not act transitively on the deformed
flag varieties even in the case of Grassmanians.
\end{rem}

\begin{rem}
We note that though $\Fl^a_{\om_d}\simeq \Fl_{\om_d}\simeq Gr(d,n)$, the actions of the
Borel groups $B\subset G$ and $B\subset G^a$ are very different.
Let us consider the case $G=SL_2$. Then $\g^a$ is spanned by three elements $e^a$, $h^a$ and $f^a$
subject to the relations
\[
[h^a,e^a]=2e^a,\ [h^a,f^a]=-2f^a,\ [e^a,f^a]=0.
\]
Let $\la$ be a dominant weight of $\msl_2$, $\la\in\bZ_{\ge 0}$. Then $V_\la^a$ is the direct sum
of one-dimensional subspaces spanned by vectors $v_l$, $l=\la, \la-2,\dots, -\la$
such that
\[
h^a v_l=l v_l,\ f^av_l=v_{l-2},\ e^a v_l=0.
\]
Therefore, the Borel subgroup $B$ acts trivially on $\Fl^a_\la\simeq \bP^1$. For instance,
there exists one point of $\bP^1$, which is fixed by the action of the whole group $G^a$.
\end{rem}

Let us now recall the Pl\"ucker relations for $\Fl_\la$ \cite{Fu} and the deformed Pl\"ucker relations for
$\Fl^a_\la$ \cite{Fe3}.

Let $1\le d_1<\dots <d_s\le n-1$ be a sequence of increasing numbers. Then for any positive
integers $a_1,\dots,a_s$ the variety $\Fl_{a_1\om_{d_1}+\dots + a_s\om_{d_s}}$ is isomorphic to
the partial flag variety
\[
\Fl(d_1,\dots,d_s)=\{V_1\hk V_2\hk\dots\hk V_s\hk \bC^n:\ \dim V_i=d_i\}.
\]
In particular, if $s=1$, then $\Fl(d)$ is the Grassmanian $Gr(d,n)$ and for $s=n-1$
$\Fl(1,\dots,n-1)$ is the variety of the complete flags. We recall that
\[
V_{\om_d}=\Lambda^d(V_{\om_1})=\Lambda^d(\bC^n)
\]
and the embedding $Gr(d,n)\hk \bP(\Lambda^d V_{\om_1})$ is defined as follows:
a subspace with a basis $w_1,\dots,w_d$ maps to $\bC w_1\wedge\dots\wedge w_d$.
For general sequence $d_1,\dots,d_s$ one has embeddings:
\[
\Fl(d_1,\dots,d_s)\hk Gr(d_1,n)\times\dots\times Gr(d_s,n)\hk
\bP(V_{\om_{d_1}})\times\dots\times \bP(V_{\om_{d_s}}).
\]
The composition of these embeddings is called the Pl\"ucker embedding.
The image is described
explicitly in terms of Pl\"ucker relations. Namely, let $v_1,\dots,v_n$ be a
basis of $\bC^n=V_{\om_1}$. Then one gets  a basis $v_J$ of $V_{\om_d}$
$v_J=v_{j_1}\wedge\dots\wedge v_{j_d}$  labeled by sequences
$J=(1\le j_1<j_2<\dots <j_d\le n)$. Let $X_J\in V_{\om_d}^*$ be the dual basis.
We denote by the same symbols the coordinates of a vector $v\in V_{\om_d}$:
$X_J=X_J(v)$.
The image of the embedding
\[
\Fl(d_1,\dots,d_s)\hk\times_{i=1}^s \bP(V_{\om_{d_i}})
\]
is defined by the Pl\"ucker relations. These relations are labeled by a pair of
numbers $p\ge q$, $p,q\in\{d_1,\dots,d_s\}$, by a number $k$, $1\le k\le q$ and by a pair
of sequences $L=(l_1,\dots,l_p)$, $J=(j_1,\dots,j_q)$, $1\le l_\al,j_\beta\le n$. The
corresponding relation is denoted by $R^k_{L,J}$ and is given by
\begin{equation}\label{PR}
R^k_{L,J}=X_LX_J - \sum_{1\le r_1 <\dots <r_k\le p} X_{L'}X_{J'},
\end{equation}
where $L'$,$J'$ are obtained from $L$, $J$ by interchanging $k$-tuples
$(l_{r_1},\dots,l_{r_k})$ and $(j_1,\dots,j_k)$ in $L$ and $J$ respectively, i.e.
\begin{gather*}
J'=(l_{r_1},\dots,l_{r_k},j_{k+1},\dots,j_q),\\
L'=(l_1,\dots,l_{r_1-1},j_1,l_{r_1+1},\dots,l_{r_2-1},j_2,\dots,l_p).
\end{gather*}
We note that for any $\sigma\in S_d$ the equality
\[
X_{j_{\sigma(1)},\dots,j_{\sigma(d)}}=(-1)^\sigma X_{j_1,\dots,j_d}
\]
is assumed in \eqref{PR}. We denote the ideal generated by all $R^k_{L,J}$ by
$I(d_1,\dots,d_s)$.

We introduce the notation
$$\Fl^a(d_1,\dots,d_s)=\Fl^a_{\om_{d_1}+\dots + \om_{d_s}},\ 1\le d_1<\dots < d_s<n.$$
\begin{dfn}
Let $I^a(d_1,\dots,d_s)$ be an ideal in the polynomial ring in variables
$X^a_{j_1,\dots,j_d}$, $d=d_1,\dots,d_s$, $1\le j_1<\dots <j_d<n$, generated by the
elements $R^{k;a}_{L,J}$ given below.
These elements are labeled by a pair of numbers $p\ge q$, $p,q\in\{d_1,\dots,d_s\}$,
by an integer $k$, $1\le k\le q$ and by sequences
$L=(l_1,\dots,l_p)$, $J=(j_1,\dots,j_q)$, which  are arbitrary subsets of the set
$\{1,\dots,n\}$. The generating elements are given by the
formulas
\begin{equation}
R^{k;a}_{L,J}=X^a_{l_1,\dots,l_p} X^a_{j_1,\dots,j_q} -
\sum_{1\le r_1<\dots <r_k\le p}
X^a_{l'_1,\dots,l'_p} X^a_{j'_1,\dots,j'_q},
\end{equation}
where the terms of $R^{k;a}_{L,J}$ are the terms of $R^{k}_{L,J}$ \eqref{PR}
(with a superscript $a$, to be precise) such that
\begin{equation}\label{new}
\{l_{r_1},\dots,l_{r_k}\}\cap \{q+1,\dots,p\}=\emptyset.
\end{equation}
\end{dfn}

\begin{rem}
The initial term $X^a_{l_1,\dots,l_p} X^a_{j_1,\dots,j_q}$ is also subject to the condition
\eqref{new}, i.e. it is not present in
$R^{k;a}_{L,J}$ if $\{j_1,\dots,j_k\}\cap \{q+1,\dots,p\}\ne\emptyset$.
\end{rem}

\begin{example}
Let $s=1$. Then $I^a(d)=I(d)$, since there are no numbers $l$ such that
$d+1\le l\le d$ and thus $R^{k;a}_{L,J}=R^{k}_{L,J}$ (up to a superscript $a$
in the notations of variables $X_J$). Hence $\Fl^a_{\om_d}\simeq \Fl_{\om_d}$.
\end{example}

The following theorem is proved in \cite{Fe3}.
\begin{thm}
The variety $\Fl^a(d_1,\dots,d_s)\hk \times_{i=1}^s \bP(\Lambda^{d_i} \bC^n)$ is
defined by the ideal $I^a(d_1,\dots,d_s)$.
\end{thm}

\begin{example}\label{1,n-1}
Let $s=2$, $d_1=1$, $d_2=n-1$. Then the classical flag variety $\Fl(1,n-1)$ is a subvariety
in $\bP^{n-1}\times \bP^{n-1}$ defined by a single relation
\[
\sum_{i=1}^n (-1)^{i-1} X_i X_{1,\dots,i-1,i+1,\dots,n}=0.
\]
The degenerate variety $\Fl(1,n-1)$ is also a subvariety in $\bP^{n-1}\times \bP^{n-1}$,
defined by  a "degenerate" relation
\[
X^a_1 X^a_{2,\dots,n} +(-1)^{n-1} X^a_nX^a_{1,\dots,n-1}=0.
\]
\end{example}

\section{Cell decomposition}\label{sec2}
In this section we describe explicitly the image of $\Fl^a_\la$ inside the product of Grassmanians
and construct the cell decomposition of the degenerate flag varieties.
We start with the case of $\la=\om_d$.

\subsection{Cell decomposition for  Grassmanians}
Recall that $\Fl^a_{\om_d}\simeq \Fl_{\om_d}\simeq Gr(d,n)$. Given an increasing tuple
$L=(l_1<\dots<l_d)$ we set
\[
p_L=\mathrm{span}(v_{l_1},\dots,v_{l_d})\in Gr(d,n).
\]
The subspace $p_L$ is $T$-invariant.
Let $k$ be a number such that $l_k\le d< l_{k+1}$.
\begin{prop}\label{Grcells}
The orbit $G^a\cdot p_L$ is an affine cell and $Gr(d,n)$ is the disjoint union of
all such cells.
\end{prop}
\begin{proof}
Recall that $G^a$ acts on $Gr(d,n)$ via the projection to the Lie group of $\g_d$ (see \eqref{gd}).
Therefore the elements of $G^a\cdot p_L$ are exactly the subspaces of $V$ having a basis
$e_1,\dots,e_d$ of the form
\begin{gather}\label{cell1}
e_j=v_{l_j}+\sum_{i=1}^{l_j-1} a_{i,j} v_i + \sum_{i=d+1}^n a_{i,j}v_i,\ j=1,\dots,k\\
\label{cell2}
e_j=v_{l_j} + \sum_{i=d+1}^{l_j-1} a_{i,j}v_i,\ j=k+1,\dots,d.
\end{gather}
Such elements in $Gr(d,n)$ obviously form an affine cell and one has a decomposition
$Gr(d,n)=\sqcup_{L} G^a\cdot p_L$.
\end{proof}

\begin{rem}\label{combine}
Formulas \eqref{cell1} and \eqref{cell2} can be combined together as follows.
Let $[k]_+=k$ if $k>0$ and $[k]_+=k+n$ if $k\le 0$. Then each element of $G^a\cdot p_L$ has
a basis $e_1,\dots,e_d$  of the form
\begin{equation}\label{+}
e_j=v_{l_j}+\sum_{i=1}^{[l_j-d]_+-1} a_{i,j} v_{[l_j-i]_+}.
\end{equation}
\end{rem}

\begin{rem}
The orbit $G^a\cdot p_L$ can be identified with a certain cell $B\cdot p_J$ in the usual cell
decomposition of $Gr(d,n)$. Namely, define $J$ as follows:
\[
J=(l_{k+1}-d,l_{k+2}-d,\dots,l_d-d,l_1-d+n,l_2-d+n,\dots,l_k-d+n).
\]
Then the map
\[
\psi:V\to V,\  \psi(v_i)=v_{[i-d]_+}, i=1,\dots,n
\]
sends $G^a\cdot p_L$ to $B\cdot p_j$ (this is clear from the explicit description
\eqref{cell1}, \eqref{cell2}).
\end{rem}

\begin{example}
Let $n=9$, $d=4$ and $L=(2,3,6,7)$ (thus $k=2$). Then the elements of $G^a\cdot p_L$ can be
identified with the following matrices (the columns of a matrix form a basis of the corresponding
subspace):
$$
\begin{pmatrix}
* & * & 0 & 0 \\
1 & 0 & 0 & 0 \\
0 & 1 & 0 & 0 \\
0 & 0 & 0 & 0 \\
* & * & * & * \\
0 & 0 & 1 & 0 \\
0 & 0 & 0 & 1 \\
* & * & 0 & 0 \\
* & * & * & *
\end{pmatrix}
$$
Here $*$ denotes arbitrary entries and hence the number of stars coincides with the dimension of the cell.
\end{example}

\subsection{Chains of subspaces}
In this section we fix the numbers $d_1,\dots,d_s$ and write $\Fl^a$ for $\Fl^a(d_1,\dots,d_s)$.
Let $v_1,\dots,v_n$ be some basis of $V\simeq \bC^n$. For $1\le i<j\le n$ we define the projections
$pr_{i+1,j}:V\to V$ by the formula
\[
pr_{i+1,j} (\sum_{l=1}^n c_lv_l)=\sum_{l=1}^i c_lv_l + \sum_{l=j+1}^n c_lv_l.
\]

The goal of this subsection is to prove the following theorem.
\begin{thm}\label{main}
The variety $\Fl^a\hk Gr(d_1,n)\times\dots\times Gr(d_s,n)$ is formed by all sequences
$V_1,\dots,V_s$, $V_l\in Gr(d_l,n)$ such that for all $1\le l<m\le s$
\begin{equation}\label{pr}
pr_{d_l+1,d_m} V_l\hk V_m.
\end{equation}
\end{thm}

\begin{rem}
It is easy to see that the set of conditions \eqref{pr} is equivalent to
the subset with $m=l+1$, i.e. to the set of conditions
\begin{equation}\label{prs}
pr_{d_l+1,d_{l+1}} V_l\hk V_{l+1},\quad l=1,\dots,s-1.
\end{equation}
\end{rem}

\begin{lem}
Let $(V_1,\dots,V_s)\in\Fl^a$. Then conditions \eqref{pr} are satisfied.
\end{lem}
\begin{proof}
Let us first look at the big cell $\bG_a^M\cdot \bC v_\la\subset \Fl^a$.
Note that the line $\bC v_\la$ is represented by the point
\[
\times_{i=1}^s\mathrm{span}(v_1,\dots,v_{d_i})\in\times_{i=1}^s Gr(d_i,n).
\]
Take an element  $g=\exp(\sum s_{i,j} f_{i,j})\in \bG_a^M\subset G^a$.
Then one has
\[
g\cdot \spa(v_1,\dots,v_d)=\spa(v_1+\sum_{j=d}^{n-1} s_{1,j}v_{j+1},\dots,
v_d+\sum_{j=d}^{n-1} s_{d,j}v_{j+1}).
\]
Therefore conditions \eqref{pr} hold for all points from the big cell of the degenerate flag varieties.
Since $\Fl^a_\la$ is the closure of the big cell, the lemma is proved.
\end{proof}

\begin{prop}
Let  $V_1,\dots,V_s$ be a set of subspaces of $V$ satisfying \eqref{pr} with $\dim V_l=d_l$.
Then $(V_1,\dots,V_s)\in \Fl^a$.
\end{prop}
\begin{proof}
We know that the image of the embedding
$$
\Fl^a\hk \times_{i=1}^s Gr(d_i,n)\hk
\times_{i=1}^s \bP(\Lambda^{d_i}V)
$$
is defined by the set of relations $R_{J,I}^{k;a}=0$. Our goal is to prove that \eqref{pr}
implies that all the relations $R_{J,I}^{k;a}$ vanish. Fix a pair $1\le l\le m\le s$.
In what follows we denote the projection
$pr_{d_l+1,d_m}$ simply by $pr$.

Let $(V_1,\dots,V_s)$ be a collection of subspaces satisfying \eqref{pr}. Fix tuples $I=(i_1,\dots,i_l)$
and $J=(j_1,\dots,j_m)$ and a number $k$. We prove that the relation $R_{J,I}^{k;a}$ vanishes on
$(V_1,\dots,V_s)$. Without loss of generality we assume that $i_1,\dots,i_k\notin [d_l+1,d_m]$.
We also rearrange the entries of $I$ in such a way that the elements from $I\cap [d_l+1,d_m]$ are concentrated
at the end of $I$, i.e. there exists a number $b$ such that
\[
i_1,\dots,i_b\notin [d_l+1,d_m],\quad i_{b+1},\dots,i_l\in [d_l+1,d_m].
\]
Obviously, $b\ge k$. Let $l-c=\dim(\ker pr\cap V_l)$.  We fix a basis $e_1,\dots,e_l$ of $V_l$ such that
$ pr e_1,\dots, pr e_c$ is a basis of $pr V_l$ and $e_{c+1},\dots,e_l$ form a basis of $\ker pr\cap V_l$.
We denote by $a_{s,t}$ the coefficients of the expansion of $e_s$ in terms of $v_t$:
\[
e_q=\sum_{r=1}^l a_{r,q} v_r.
\]
The idea of the proof is to use the following decomposition of a Pl\"ucker coordinate $X_I$:
\begin{equation}\label{decomp}
X_I=\sum_{1\le \al_1<\dots <\al_{l-b}\le l} \pm a_{i_{b+1},\al_1}\dots a_{i_{l},\al_{l-b}}X_{i_1,\dots,i_b}.
\end{equation}
Here $X_{i_1,\dots,i_b}$ is the $(i_1,\dots,i_b)$-th Pl\"ucker coordinate of the vector space
$\mathrm{span}(e_{\beta_1},\dots,e_{\beta_b})$, where the set of $\beta$'s is complementary to the
set of $\al$'s, i.e.
\[
\{\beta_1,\dots,\beta_{b}\}\cup \{\al_1,\dots,\al_{l-b}\}=\{i_1,\dots,i_l\}.
\]
The decomposition \eqref{decomp} induces the decomposition of the relation $R^{k;a}_{J,I}$, such that each term
can be shown to vanish. Note that if $b>c$ then $X_I$ vanishes on $V_l$. We thus assume that $b\le c$.

Define the subspace
\[
E_\beta=pr (\mathrm{span} (e_{\beta_1},\dots, e_{\beta_b})).
\]
We know that $E_\beta\hk V_m$. In addition, the coordinates $X_{(i_1,\dots,i_b)}$ of
the space $\mathrm{span} (e_{\beta_1},\dots, e_{\beta_b})$ coincide with the Pl\"ucker coordinates
$Y_{(i_1,\dots,i_b)}$ of $E_\beta$, because $i_1,\dots,i_b\notin [d_l+1,d_m]$ (we are using the notations
$Y_I$ to distinguish between Pl\"ucker coordinated of different spaces). Since $E_\beta\hk V_m$,
the classical relations $R^k_{J,(i_1,\dots,i_b)}$ vanish on the pair $(E_\beta,V_m)$. Since
\[
E_\beta\hk\mathrm{span}(v_1,\dots,v_{d_l},v_{d_m+1},\dots,v_n),
\]
a Pl\"ucker coordinate $Y_{q_1,\dots,q_b}$ of $E_\beta$ vanishes unless non of the indices $q_\bullet$ are
between $d_l+1$ and $d_m$. Hence the degenerate Pl\"ucker relation $R^{k:a}_{J,(i_1,\dots,i_b)}$ also vanishes on
$(E_\beta,V_m)$. Note also that the decomposition \eqref{decomp} induces the decomposition
\[
R^{k;a}_{J,I}=\sum_{1\le \al_1<\dots <\al_{l-b}\le l} \pm a_{i_{b+1},\al_1}\dots a_{i_{l},\al_{l-b}}
R^{k;a}_{J,(i_{\beta_1},\dots,i_{\beta_b})}.
\]
But as we have shown above, each of the relations $R^{k;a}_{J,(i_{\beta_1},\dots,i_{\beta_b})}$ vanishes
on $(V_l,V_m)$. Hence so does $R^{k;a}_{J,I}$.
\end{proof}

\begin{example}
Let $\la=\om_1+\om_{n-1}$, i.e. $s=2$, $d_1=1$, $d_2=n-1$. Then the image of $\Fl^a(1,n-1)$ inside
$Gr(1,n)\times Gr(n-1,n)$ is formed by all pairs $V_1,V_2$ such that $pr_{2,n-1} V_1\hk V_2$.
Since $pr_{2,n-1} V_1\hk\spa(v_1,v_n)$, the image of the embedding
$\Fl^a(1,n-1)\hk \bP^{n-1}\times \bP^{n-1}$is defined by a single relation
\[
X^a_1X^a_{2,\dots,n} + (-1)^{n-1} X_n^a X^a_{1,\dots,n-1}=0,
\]
which agrees with Example \ref{1,n-1}.
 \end{example}

\begin{cor}
Theorem \ref{main} is true.
\end{cor}

\begin{cor}
Let $I^1,\dots,I^s$, $I^l\subset \{1,\dots,n\}$ be a collection of tuples
such that the cardinality of $I^l$ is $d_l$. Then a point $p_{I^1}\times\dots\times p_{I^s}$ belongs to $\Fl^a$
if and only if
\begin{equation}\label{ngen}
I^l\setminus \{d_l+1,\dots,d_{l+1}\}\subset I^{l+1}.
\end{equation}
\end{cor}

\begin{example}\label{complete}
Consider the case of the complete flags: $s=n-1$, $d_l=l$. Set $pr_l=pr_{l,l}$. Then  the embedding of $\Fl^a$
into the product of Grassmanians is defined by the conditions
\begin{equation}
pr_{l+1} V_l\hk V_{l+1},\ l=1,\dots,n-2
\end{equation}
and the conditions \eqref{ngen} read as
$I^l\setminus \{l+1\}\subset I^{l+1}$ for $l=1,\dots,n-2$.
\end{example}

\subsection{Cells for $\Fl^a$}
Recall that the cell decomposition for a Grassmanian is given by the $G^a$-orbits
of the torus fixed points. However this is not true for the case of general $\Fl^a_\la$.
Moreover, the number of $G^a$-orbits can be infinite.
The simplest example is as follows.
\begin{example}
Let $n=4$, $\la=\om_1+\om_3$. Then $\Fl^a_\la$ is embedded into $\bP^3\times \bP^3$
(two Grassmanians for $\msl_4$) with the coordinates
$(x_1:x_2:x_3:x_4)$ and $(x_{123}:x_{124}:x_{133}:x_{234})$. The variety $\Fl^a_{\om_1+\om_3}$
is defined by a single relation
$x_1x_{234}-x_4x_{123}=0$. Therefore, $\Fl^a_{\om_1+\om_3}$ contains the product
$\bP^2\times \bP^2$ defined by $x_1=x_{123}=0$. We note that the subgroup $\bG_a^{6}$ of $G^a$
acts trivially on this $\bP^2\times \bP^2$ (the PBW-degree in both $V_{\om_1}$ and $V_{\om_3}$
is at most one). Therefore, we are left with an action of the Borel subgroup.
Let $w_1,w_2,w_3,w_4$ and $w_{123},w_{124},w_{134},w_{234}$ be the standard bases for $V_{\om_1}$
and $V_{\om_3}$.
The group $B$ acts on the span of $w_2,w_3,w_4$ (resp. on the span of $w_{124},w_{134},w_{234}$)
as on the quotient of the vector representation (resp. the dual vector representation)
by $\bC w_1$ (resp. $\bC w_{123}$). It is easy to see that the corresponding
$B$-action on $\bP^2\times \bP^2$ has infinitely many orbits.
\end{example}

In the following proposition we describe the cell decomposition for $\Fl^a=\Fl^a(d_1,\dots,d_s)$.
\begin{prop}
Let ${\bf I}=(I^1,\dots,I^s)$ be a set of sequences satisfying the condition \eqref{ngen}.
Then there exists a cell decomposition $\Fl^a=\sqcup_{{\bf I}} C_{{\bf I}}$, where
\[
C_{{\bf I}}=(G^a\cdot p_{I^1}\times \dots \times G^a\cdot p_{I^s})\cap \Fl^a.
\]
In other words, a cell is given by the intersection of the degenerate flag variety, embedded into
the product of Grassmanians, with the product of the corresponding cells in $Gr(d_i,n)$.
\end{prop}
\begin{proof}
In Theorem \ref{Poincare} we compute the dimensions of $C_{{\bf I}}$. In the proof we construct
explicitly the coordinates on $C_{{\bf I}}$ thus showing that $C_{{\bf I}}$ is a cell.
\end{proof}

\section{The median Genocchi numbers}\label{sec3}
\subsection{Combinatorics}
Let $h_n$ be the normalized Genocchi numbers of the second kind. They are also referred to as the
normalized median Genocchi numbers.  These numbers have several definitions
(see \cite{De}, \cite{Du}, \cite{Kr}, \cite{S}).
The first several $h_n$'s are as follows: $1,2,7,38,295,3098$.
We first briefly recall definitions of these numbers.

We start with the Dellac definition (see \cite{De}). Consider a rectangle with $n$ columns and
$2n$ rows. It contains $n\times 2n$ boxes labeled by pairs $(l,j)$, where $l=1,\dots,n$ is the number
of a column and $j=1,\dots,2n$ is the number of  a row. A Dellac  configuration $D$ is a subset of boxes,
subject to the following conditions:
\begin{itemize}
\item each column contains exactly two boxes from $D$,
\item each row contains exactly one box from $D$,
\item if the $(l,j)$-th box is in $D$, then $l\le j\le n+l$.
\end{itemize}
Let $DC_n$ be the set of such configurations. Then the number of elements in $DC_n$ is equal to $h_n$.

We list all Dellac's configurations for $n=3$.
We specify boxes in a configuration by putting  fat dots inside.

\begin{equation}\label{n=3}
\begin{picture}(30,60)
\put(0,0){\line(1,0){30}}
\put(0,10){\line(1,0){30}}
\put(0,20){\line(1,0){30}}
\put(0,30){\line(1,0){30}}
\put(0,40){\line(1,0){30}}
\put(0,50){\line(1,0){30}}
\put(0,60){\line(1,0){30}}

\put(0,0){\line(0,1){60}}
\put(10,0){\line(0,1){60}}
\put(20,0){\line(0,1){60}}
\put(30,0){\line(0,1){60}}

\put(2,2){$\bullet$}
\put(2,12){$\bullet$}
\put(12,22){$\bullet$}
\put(12,32){$\bullet$}
\put(22,42){$\bullet$}
\put(22,52){$\bullet$}
\end{picture} \quad
\begin{picture}(30,60)
\put(0,0){\line(1,0){30}}
\put(0,10){\line(1,0){30}}
\put(0,20){\line(1,0){30}}
\put(0,30){\line(1,0){30}}
\put(0,40){\line(1,0){30}}
\put(0,50){\line(1,0){30}}
\put(0,60){\line(1,0){30}}

\put(0,0){\line(0,1){60}}
\put(10,0){\line(0,1){60}}
\put(20,0){\line(0,1){60}}
\put(30,0){\line(0,1){60}}

\put(2,2){$\bullet$}
\put(2,12){$\bullet$}
\put(12,22){$\bullet$}
\put(12,42){$\bullet$}
\put(22,32){$\bullet$}
\put(22,52){$\bullet$}
\end{picture} \quad
\begin{picture}(30,60)
\put(0,0){\line(1,0){30}}
\put(0,10){\line(1,0){30}}
\put(0,20){\line(1,0){30}}
\put(0,30){\line(1,0){30}}
\put(0,40){\line(1,0){30}}
\put(0,50){\line(1,0){30}}
\put(0,60){\line(1,0){30}}

\put(0,0){\line(0,1){60}}
\put(10,0){\line(0,1){60}}
\put(20,0){\line(0,1){60}}
\put(30,0){\line(0,1){60}}

\put(2,2){$\bullet$}
\put(2,12){$\bullet$}
\put(12,32){$\bullet$}
\put(12,42){$\bullet$}
\put(22,22){$\bullet$}
\put(22,52){$\bullet$}
\end{picture} \quad
\begin{picture}(30,60)
\put(0,0){\line(1,0){30}}
\put(0,10){\line(1,0){30}}
\put(0,20){\line(1,0){30}}
\put(0,30){\line(1,0){30}}
\put(0,40){\line(1,0){30}}
\put(0,50){\line(1,0){30}}
\put(0,60){\line(1,0){30}}

\put(0,0){\line(0,1){60}}
\put(10,0){\line(0,1){60}}
\put(20,0){\line(0,1){60}}
\put(30,0){\line(0,1){60}}

\put(2,2){$\bullet$}
\put(2,22){$\bullet$}
\put(12,12){$\bullet$}
\put(12,32){$\bullet$}
\put(22,42){$\bullet$}
\put(22,52){$\bullet$}
\end{picture}\quad
\begin{picture}(30,60)
\put(0,0){\line(1,0){30}}
\put(0,10){\line(1,0){30}}
\put(0,20){\line(1,0){30}}
\put(0,30){\line(1,0){30}}
\put(0,40){\line(1,0){30}}
\put(0,50){\line(1,0){30}}
\put(0,60){\line(1,0){30}}

\put(0,0){\line(0,1){60}}
\put(10,0){\line(0,1){60}}
\put(20,0){\line(0,1){60}}
\put(30,0){\line(0,1){60}}

\put(2,2){$\bullet$}
\put(2,22){$\bullet$}
\put(12,12){$\bullet$}
\put(12,42){$\bullet$}
\put(22,32){$\bullet$}
\put(22,52){$\bullet$}
\end{picture}\quad
\begin{picture}(30,60)
\put(0,0){\line(1,0){30}}
\put(0,10){\line(1,0){30}}
\put(0,20){\line(1,0){30}}
\put(0,30){\line(1,0){30}}
\put(0,40){\line(1,0){30}}
\put(0,50){\line(1,0){30}}
\put(0,60){\line(1,0){30}}

\put(0,0){\line(0,1){60}}
\put(10,0){\line(0,1){60}}
\put(20,0){\line(0,1){60}}
\put(30,0){\line(0,1){60}}

\put(2,2){$\bullet$}
\put(2,32){$\bullet$}
\put(12,12){$\bullet$}
\put(12,22){$\bullet$}
\put(22,42){$\bullet$}
\put(22,52){$\bullet$}
\end{picture}\quad
\begin{picture}(30,60)
\put(0,0){\line(1,0){30}}
\put(0,10){\line(1,0){30}}
\put(0,20){\line(1,0){30}}
\put(0,30){\line(1,0){30}}
\put(0,40){\line(1,0){30}}
\put(0,50){\line(1,0){30}}
\put(0,60){\line(1,0){30}}

\put(0,0){\line(0,1){60}}
\put(10,0){\line(0,1){60}}
\put(20,0){\line(0,1){60}}
\put(30,0){\line(0,1){60}}

\put(2,2){$\bullet$}
\put(2,32){$\bullet$}
\put(12,12){$\bullet$}
\put(12,42){$\bullet$}
\put(22,22){$\bullet$}
\put(22,52){$\bullet$}
\end{picture}\ .
\end{equation}

The Dellac definition is the earliest one, but the most well-known definition is via the Seidel triangle.
The Seidel triangle is of the form
$$
\begin{matrix}
1 & & & & \\
1 & & & & \\
1 & 1 & & & \\
2 & 1 & & & \\
2 & 3 & 3 & & \\
8 & 6 & 3 & & \\
8 & 14 & 17 & 17 & \\
56 & 48 & 34 & 17 & \\
56 & 104 & 138 & 155 & 155
\end{matrix}
$$
By definition, the triangle is formed by the numbers $G_{k,n}$
($n$ is the number of a row and $k$ is the number of a column)
with $n=1,2,\dots$ and $1\le k\le \frac{n+1}{2}$,
subject to the relations $G_{1,1}=1$ and
\[
G_{k,2n}=\sum_{i\ge k} G_{i,2n-1},\ G_{k,2n+1}=\sum_{i\le k} G_{i,2n}.
\]
The numbers $G_{n,2n-1}$ are called the Genocchi numbers of the first kind and the numbers
$G_{1,2n}$ are called the Genocchi numbers of the second kind (or the median Genocchi numbers).
Barsky \cite{Ba} and then Dumont \cite{Du} proved that the number $G_{1,2n+2}$ is divisible by
$2^n$. The normalized median Genocchi numbers $h_n$ are defined as the corresponding ratios:
$h_n=G_{1,2n+2}/2^n$.

In \cite{Kr} Kreweras suggested another description of the numbers $h_n$. Namely, a permutation
$\sigma\in S_{2n+2}$ is called a normalized Dumont permutation of the second kind if the following conditions
are satisfied:
\begin{itemize}
 \item $\sigma(k)<k$ if $k$ is even,
\item $\sigma(k)>k$ if $k$ is odd,
\item $\sigma^{-1}(2k)<\sigma^{-1}(2k+1)$ for $k=1,\dots,n$.
\end{itemize}
The set of such permutations is denoted by $PD2N_n$ (P for permutations, D for Dumont, $2$ for the second
kind and N for normalized). According to Kreweras, the number of elements
of $PD2N_n$ is equal to $h_n$. In Proposition \ref{Dellac} we show that the definitions of Dellac and
Kreweras are equivalent (this seems to be known to expert -- see \cite{G}, \cite{S}, but we were not
able to find a proof in the literature).

In the following proposition we show that the conditions from Example \ref{complete}
give rise to a new definition of the numbers $h_n$.
\begin{prop}\label{fullgen}
The number of tuples $I^1,\dots, I^{n-1}$, with $I^l\subset \{1,\dots,n\}$, $\# I^l=l$ subject to the condition
\begin{equation}\label{combpr}
I^{l-1}\setminus\{l\}\subset I^l,\ l=2,\dots,n-1
\end{equation}
is equal to $h_n$.
\end{prop}
\begin{proof}
Let $\bar h_n$ be the number of tuples as above. We compare $\bar h_n$ with the Dellac definition
of $h_n$. Given a set $I^1,\dots,I^{n-1}$ subject to the condition \eqref{combpr}, we construct
the corresponding Dellac's configuration $D$ and then prove that this map is one-to-one.
The rule is as follows. Let us explain what are the boxes of $D$ in the $l$-th column.

First, suppose $l\notin I^{l-1}$. Then because of the condition \eqref{combpr} the difference
$I^l\setminus I^{l-1}$ contains exactly one number $j$.  There are two cases:
\begin{itemize}
\item If $j>l$, then $D$ contains boxes $(l,l)$ and $(l,j)$.
\item If $j\le l$, then $D$ contains boxes $(l,l)$ and $(l,j+n)$.
\end{itemize}
Now, suppose $l\in I^{l-1}$.  Then either $l\in I^l$, or $L\notin I^l$.
If $l\in I^l$, then $I^l\setminus I^{l-1}$ contains exactly one number $j$.
There are two cases:
\begin{itemize}
\item If $j>l$, then $D$ contains boxes $(l,l+n)$ and $(l,j)$.
\item If $j\le l$, then $D$ contains boxes $(l,l+n)$ and $(l,j+n)$.
\end{itemize}
Finally, let $l\in I^{l-1}$ and $l\notin I^l$. Then $I^l\setminus I^{l-1}$ contains exactly
two numbers $j_1$ and $j_2$. There are four variants:
\begin{itemize}
\item If $j_1>l$ and $j_2>l$, then $D$ contains boxes $(l,j_1)$ and $(l,j_2)$.
\item If $j_1>l$ and $j_2\le l$, then $D$ contains boxes $(l,j_1)$ and $(l,n+j_2)$.
\item If $j_1\le l$ and $j_2>l$, then $D$ contains boxes $(l,j_1+n)$ and $(l,j_2)$.
\item If $j_1\le l$ and $j_2\le l$, then $D$ contains boxes $(l,j_1+n)$ and $(l,j_2+n)$.
\end{itemize}
This rule explains how to pick boxes in columns from $1$ to $n-1$. To complete the configuration
we simply pick two boxes in the last column in the unique way to make $D$ a configuration.

In order to prove that this map is a bijection, we construct the inverse map.
Let $D$ be a Dellac configuration. We define $I^l$ inductively. First, let $l=1$.
Then the box $(1,1)$ necessarily belongs to $D$. Let $j>1$ and $D$ contains $(1,j)$. Then
if $j=n+1$, then $I^1=(1)$. Otherwise $I^1=(j)$.

Now assume that $I^{l-1}$ is already defined.
First, suppose that the $(l,l)$-th box belongs to $D$. Then there exists one more box
$(l,j)$ in $D$ with $n+l\ge j>l$. If $j\le n$ we set $I^l=I^{l-1}\cup \{j\}$.
Otherwise, we set $I^l=I^{l-1}\cup \{j-n\}$.
Second, suppose that the $(l,l)$-th box does not belong to $D$.
Since the $l$-th row of $D$ contains exactly one box, there exists $l_1<l$ such that the $(l_1,l)$-th box  belongs to $D$.
Therefore, $l\subset I^{l-1}$.
There exist exactly two boxes $(l,j_1)$ and $(l,j_2)$ in $D$ in the $l$-th column.
Then we set $I^l=I^{l-1}\setminus \{l\}\cup \{\bar j_1, \bar j_2\}$, where $\bar j=j$, if $j\le n$ and
$\bar j=j-n$ otherwise.
\end{proof}

\begin{example}
Let $n=3$.
The pairs $I^1,I^2$,  corresponding to the Dellac configurations \eqref{n=3} are as
follows (the order is the same as on picture \eqref{n=3}):
\begin{gather*}
\{(2), (13)\},\quad \{(2), (23)\},\quad \{(2), (12)\},\quad \{(3), (13)\},\\
\{(3), (23)\},\quad \{(1), (13)\},\quad \{(1), (12)\}.
\end{gather*}
\end{example}

We now compare the definitions by Dellac and by Kreweras.
\begin{prop}\label{Dellac}
The number of elements in $PD2N_n$ is equal to the number of elements in $DC_n$.
\end{prop}
\begin{proof}
We construct a bijection $A:PD2N_n\to DC_n$. Let $\sigma\in PD2N_n$.
We determine the boxes in the $k$-th column of $A(\sigma)$ using the values of
$\sigma^{-1} (2k)$ and $\sigma^{-1} (2k+1)$.

Let us start with $k=1$. We note that
$\sigma (2)=1$, $\sigma(4)$ is equal to $2$ or to $3$.
In addition, $\sigma^{-1}(2)=1$ or $4$ and the possible values of $\sigma^{-1}(3)$ are
$4,6,\dots,2n+2$.
Therefore, all possible values of the pair $(\sigma^{-1}(2), \sigma^{-1}(3))$ are
as follows:
\[
(1,4),\ (4,6), \ (4,8), \dots, (4,2n+2).
\]
If the first possibility occurs, then by definition the first column of $A(\sigma)$  contains boxes
$(1,1)$ (as any Dellac's configuration) and $(1,n+1)$. If
$\sigma^{-1}(2)=4$ and $\sigma^{-1}(3)=2l+2$, then the first column of $A(\sigma)$  contains boxes
$(1,1)$ and $(1,l)$.

Now let us consider the case $k=n$.
We note that
$\sigma (2n+1)=2n$, $\sigma(2n-1)$ is equal to $2n$ or to $2n+1$.
In addition, $\sigma^{-1}(2n+1)=2n+2$ or $2n-1$ and the possible values of
$\sigma^{-1}(2n)$ are $1,3,\dots,2n-1$.
Therefore, all possible values of the pair $(\sigma^{-1}(2n), \sigma^{-1}(2n+1))$ are
as follows:
\[
(2n-1,2n+2),\ (1,2n-1), \ (3,2n-1), \dots, (2n-3,2n-1).
\]
If the first possibility occurs, then by definition the $n$-th column of $A(\sigma)$  contains boxes
$(n,2n)$ (as any Dellac's configuration) and $(n,n)$. If
$$(\sigma^{-1}(2n), \sigma^{-1}(2n+1))=(2l-1,2n-1),$$
then the first column of $A(\sigma)$  contains boxes
$(n,2n)$ and $(n,n+l)$.

Finally,  take $k=2,\dots,n-1$. We note that the possible values of $\sigma^{-1}(2k)$ are
$1,3,\dots,2k-1,2k+2,\dots,2n$. Also, the possible values of $\sigma^{-1}(2k+1)$ are
$3,5,\dots,2k-1,2k+2,\dots,2n, 2n+2$.
We now define the $k$-th column of $A(\sigma)$ as follows:
\begin{enumerate}
\item \label{i} If the pair $(\sigma^{-1}(2k),\sigma^{-1}(2k+1))$ contains $2l-1$, $l=1,\dots,k$,
then the $k$-th column of $A(\sigma)$ contains a box $(k,n+l)$.
\item \label{ii} If the pair $(\sigma^{-1}(2k),\sigma^{-1}(2k+1))$ contains $2l+2$, $l=k,\dots,n$,
then the $k$-th column of $A(\sigma)$ contains a box $(k,l)$.
\end{enumerate}

We note that  $A(\sigma)\in DC_n$. In fact, by definition any column of
$A(\sigma)$ contains exactly two boxes and every row contains exactly
one box (this follows from the definition above and because $\sigma$ is one-to-one).
In order to prove that $A$ is a bijection it suffices to note that formulas \eqref{i} and
\eqref{ii} allow to construct explicitly the map $A^{-1}$.
\end{proof}

\begin{example}
Let $n=3$. The elements of $PD2N_3$ corresponding to the Dellac configurations on picture
\eqref{n=3} are as follows (the order is the same as on picture \eqref{n=3}):
\begin{gather*}
(41627385), \quad (61427385),\quad (41526387),\quad (41627583),\\
(61427583),\quad (21637485),\quad (21436587).
\end{gather*}
\end{example}

We recall that the main ingredient for the Kreweras construction
of $PD2N_n$ is the following triangle:
\[
\begin{picture}(120,85)
\put(0,0){295}
\put(20,0){552}
\put(40,0){702}
\put(60,0){702}
\put(80,0){552}
\put(100,0){295}
\put(10,15){38}
\put(30,15){69}
\put(50,15){81}
\put(70,15){69}
\put(90,15){38}
\put(20,30){7}
\put(40,30){12}
\put(60,30){12}
\put(80,30){7}
\put(30,45){2}
\put(50,45){3}
\put(70,45){2}
\put(40,60){1}
\put(60,60){1}
\put(50,75){1}
\end{picture}
\]
The rule is as follows: denote the numbers in the $n$-th line by $h_{n,1},\dots,h_{n,n}$.
For example, $h_{4,2}=12$.
Then the Kreweras triangle is defined by
\begin{gather*}
h_{n,1}=h_{n-1,1}+\dots+ h_{n-1,n-1},\ h_{n,2}=2h_{n,1} - h_{n-1,1},\\
h_{n,k}=2h_{n,k-1}-h_{n,k-2}-h_{n-1,k-2}-h_{n-1,k-1},\ k\ge 3.
\end{gather*}
Kreweras proved that $h_{n+1,1}$ is the $n$-th Genocchi number $h_n$ and in general $h_{n+1,k}$
is the number of the normalized Dumont permutations $\sigma\in S_{2n+2}$ of the second kind such that
$\sigma(1)=2k$.
The following is an immediate corollary from the explicit bijections above.
\begin{cor}
The number of the Dellac configurations $D\in DC_n$ such that
$\min\{i:\ (i,n+1)\in D\}=k$ is equal to $h_{n,k}$.
The number of tuples $I^1,\dots,I^{n-1}$ subject to the condition $I^{l-1}\setminus\{l\}\subset I^l$
with an extra condition  $\min\{j:\ 1\in I^j\}=k$ is equal to $h_{n,k}$.
\end{cor}

\subsection{The Poincar{\' e} polynomials.}\label{Pp}
For a tuple ${\bf I}=(I^1,\dots,I^{n-1})$ subject to the relation $I^{l-1}\setminus \{l\}\subset I^l$
we denote by $D_{\bf I}$ the corresponding Dellac configuration. For a Dellac configuration
$D\in DC_n$ we define the length $l(D)$ of $D$ as the number of pairs $(l_1,j_1)$, $(l_2,j_2)$ such that
the boxes $(l_1,j_1)$ and $(l_2,j_2)$ are both in $D$ and $l_1<l_2$, $j_1>j_2$.
We call such a pair of boxes $(l_1,j_1)$,  $(l_2,j_2)$ a disorder.
This definition resembles the definition of the length of a permutation.
We note that in the classical case the dimension of a cell attached to  a permutation $\sigma$ in a
flag variety is equal to the number of pairs $j_1<j_2$ such that $\sigma (j_1)>\sigma(j_2)$
(which equals to the length of $\sigma$).

\begin{thm}\label{Poincare}
The dimension of a cell $C_{\bf I}$ is equal to $l(D_{\bf I})$.
\end{thm}
\begin{proof}
We prove the dimension formula by constructing explicitly the coordinates on the cell $C_{\bf I}$. Let
\[
{\bf I}=(I^1,\dots,I^{n-1}),\
I^d=(i^d_1< \dots < i^d_d).
\]
Recall the description of the cells $C_{I^d}\subset Gr(d,n)$ from Proposition \ref{Grcells}.
Using this description we construct the coordinates on $C_{\bf I}$
inductively on $d$. Let $(V_1,\dots,V_{n-1})\in C_{\bf I}$.
For a number $k$ we set $[k]_+=k$ if $k>0$ and $[k]_+=k+n$ if $k\le 0$.

We start with $d=1$. An element $V_1\in C_{I^1}$ is
of the form $\bC e^1_1$ with
\[
e^1_1=v_{i^1_1}+ a^1_1 v_{[i^1_1-1]_+} + \dots + a^1_{[i^1_1-1]_+ -1} v_{2}
\]
(see Remark \ref{combine}).
We state that $[i^1_1-1]_+ -1$ (which is exactly the number of the degrees of freedom we have so far)
is exactly the number of boxes $(l,j)\in D_{\bf I}$ such that
$l>1$ and $j<i_1^1$ (note that the box $(1,1)$ is necessarily in $D_{\bf I}$, but it does not add
anything to the length of $D_{\bf I}$, since for any $(l.j)\in D_{\bf I}$ with $l>1$ we have $j>1$).
In fact, the first column of $D_{\bf I}$ contains boxes in the first row and in the
$([i^1_1-1]_+ +1)$-st row (see the proof of Proposition \ref{fullgen}).
Since any row of $D_{\bf I}$ contains exactly one box, the rows number $2,\dots, [i^1_1-1]_+$ are
occupied by boxes in the columns from $2$ to $n$. Therefore, the box $(1,[i^1_1-1]_+ +1)$ produces
exactly  $[i^1_1-1]_+ -1$ disorders.

The second step is to construct the coordinates on those subspaces from
$C_{I^2}$ which contain $pr_2 V_1$. There are two possibilities: either $i^1_1=2$ or
$i^1_1\ne 2$. In the first case  the condition $pr_2 V_1\hk V_2$ is empty.
Therefore, we have to choose two basis vectors $e^2_1,e^2_2$ of $V_2\in C_{I^2}$, with the coordinates
\begin{gather*}
e^2_1=v_{i^2_1}+ a^1_1 v_{[i^2_1-1]_+} + \dots + a^1_{[i^2_1-2]_+ -1} v_3,\\
e^2_2=v_{i^2_2}+ a^2_1 v_{[i^2_2-1]_+} + \dots + a^2_{[i^2_2-2]_+ -2} v_3.
\end{gather*}
We note that the number of coefficients of $e^2_2$ is $[i^2_2-2]_+ -2$, because
$i_1^2<i_2^2$ and hence adding appropriately normalized vector $e_1^2$ one can vanish the coefficient
of $e^2_2$ in front of $v_{i^2_1}$. We note that since $i^1_1=2$, the second column of $D_{\bf I}$
contains boxes in the rows $([i^2_1-2]_+ +2)$ and $([i^2_1-2]_+ +2)$
(see the proof of Proposition \ref{fullgen}). We state that
$[i^2_1-2]_+ -1 + [i^2_2-2]_+ -2$ (the number of degrees of freedom we have fixing the vectors $e^2_1$ and
$e^2_2$)
is exactly the number of boxes in the columns $3,4,\dots, n$, having
disorders with boxes in the second column. In fact, each row from $3$ to $[i^2_1-2]_+ -1$ contains
one box in the columns $3$ and greater (recall $i^1_1=2$). This produces $[i^2_1-2]_+ -1$ disorders
with the box $(2, [i^2_1-2]_+ - 1)$. Similarly, we obtain $[i^2_2-2]_+ -2$ disorders with the second box in the
second column.

Now assume $i_1^1\ne 2$. Then the space $pr_2 V_1$ is nontrivial and spanned by a single
vector $e^2_1=pr_2 e^1_1$. Therefore in order to specify $V_2$ we need to fix one more vector $e_2^2$ such
that $\mathrm{span}(e_1^2,e_2^2)\in C_{I^2}$. Recall that since $i_1^1\ne 2$
we have $I^2\setminus I^1=\{j\}$. Also, the second column
of $D_{\bf I}$ contains boxes in the second row and  in the row number $[j-2]_+ +2$
(see the proof of Proposition \ref{fullgen}). The box
$(2,2)$ does not produce any disorder with boxes in the columns greater than $2$. As for the box
$(2,[j-2]_+ +2)$, the number of disorders it produces is equal to the number of degrees of freedom
of choosing the vector $e^2_2$ (the argument is very similar to the ones above in the case $i^1_1=2$).

Now let us consider the general induction step. Assume that we have already computed
the number of degrees of freedom while fixing the subspaces $V_1,\dots, V_{d-1}$. Our goal is to show
that the number of degrees of freedom of $V_d$ is equal to the number of disorders produced
by the boxes in the $d$'th column with the boxes in columns $l$ with $l>d$. As in the previous case, one
has to consider two cases: $d\in I^{d-1}$ and $d\notin I^{d-1}$. The proof is very similar to the
one in the case $d=2$ and we omit it.
\end{proof}

\begin{cor}
The Poincar\' e polynomial $P_n(t)=P_{\Fl^a}(t)$ is given by
\[
P_n(t)=\sum_{D\in DC_n} t^{2l(D)}.
\]
\end{cor}

Let $q=t^2$. Then $P_n$ are polynomials in $q$ with $P_n(1)=h_n$.
Thus the Poincar\' e polynomials of the degenerate flag varieties provide a natural
$q$-version of the normalized median Genocchi numbers (it would be interesting to compare our
$q$-version with the one in \cite{HZ}).
\begin{example}
The first four polynomials $P_n(q)$ are as follows:
\begin{gather*}
P_1(q)=1, \qquad P_2(q)=1+q,\\
P_3(q)=1+ 2q+ 3q^2+q^3,\\
P_4(q)=1+3q+7q^2+10q^3+10q^4+6q^{5}+q^{6}.
\end{gather*}
\end{example}

\section*{Acknowledgments}
This work was partially supported
by the Russian President Grant MK-281.2009.1,  RFBR Grant 09-01-00058,
by grant Scientific Schools 6501.2010.2 and by the Dynasty Foundation.

\end{document}